\documentclass[11pt]{amsart}
\usepackage{amssymb, amsthm, graphicx, epsfig, graphs, verbatim}
\newtheorem{thm}{Theorem}[section]
\newtheorem{cor}[thm]{Corollary}
\newtheorem{lem}[thm]{Lemma}
\newtheorem{prop}[thm]{Proposition}

\theoremstyle{definition}
\newtheorem{defn}[thm]{Definition}

\newtheorem{rem}[thm]{Remark}

\numberwithin{equation}{section}

\newcommand{\bfz}{{\mathbb {Z}}}

\newcommand{\s}{\mathbf s}

\renewcommand{\t}{\mathbf t}

\newcommand{\Z}{\mathbb Z}

\newcommand{\Q}{\mathbb Q}

\newcommand{\del}{\partial}

\newcommand{\co}{\colon\thinspace} 
\newcommand{\hf}{{{\widehat {HF}}}} 
 
\newcommand{\Li}{\mathbb {L}}

\newcommand{\mk}{{\overline {K}}}
\newcommand{\ml}{{\overline {L}}}
\newcommand{\T}{{\mathcal {T}}^+}
\newcommand{\chat}{{\hat {c}}}

\DeclareMathOperator{\tb}{tb}

\DeclareMathOperator{\rk}{rk}

\begin{document}
\mathsurround=1pt 
\title{Notes on the contact Ozsv\'ath--Szab\'o invariants}

\author{Paolo Lisca}
\address{Dipartimento di Matematica\\
Universit\`a di Pisa \\I-56127 Pisa, ITALY} 
\email{lisca@dm.unipi.it}

\author{Andr\'{a}s I. Stipsicz}
\address{R\'enyi Institute of Mathematics\\
Hungarian Academy of Sciences\\
H-1053 Budapest\\ 
Re\'altanoda utca 13--15, Hungary and 
Institute for Advanced Study, Princeton, NJ}
\email{stipsicz@math-inst.hu and stipsicz@math.ias.edu}

\begin{abstract}
In this paper we prove various results on contact structures obtained
by contact surgery on a single Legendrian knot in the standard contact
three--sphere. Our main tool are the contact Ozsv\'ath--Szab\'o
invariants.
\end{abstract}

\thanks{The first author was partially supported by MURST, and he is a
  member of EDGE, Research Training Network HPRN-CT-2000-00101,
  supported by The European Human Potential Programme. The authors
  would like to thank Peter Ozsv\'ath and Zolt\'an Szab\'o for many
  useful discussions regarding their joint work. The second author was
  partially supported OTKA T037735.}

\maketitle

\section{Introduction}\label{s:intro}
According to a recent result of Ding and Geiges \cite{DG2} any closed
contact 3--manifold is obtained by contact surgery along a Legendrian
link $\Li$ in the standard contact 3--sphere $(S^3, \xi_{st})$, where
the surgery coefficients on the individual components of $\Li$ can be
chosen to be $\pm 1$ relative to the contact framing. (For additional
discussion on this theorem see \cite{DGS}.)  It is an intriguing
question how to establish interesting properties of a contact
structure from one of its surgery presentations. More precisely, we
would like to find a way to determine whether the result of a certain
contact surgery is tight or fillable. Recall that contact
$(-1)$--surgery (also called \emph{Legendrian surgery}) on a 
Legendrian link $\Li$ produces a Stein fillable, hence tight contact
3--manifold.

Given a Legendrian knot $K\subset (S^3, \xi_{st})$, we shall denote
the result of contact $(+1)$--surgery along $K$ by $(Y_K, \xi_K)$. A
first result, which has an elementary proof, is the following.

\begin{thm}\label{t:ot}
Let $K$ be a Legendrian knot in the standard contact three--sphere.
Assume that, for some orientation of $K$, a front projection of $K$
contains the configuration of Figure~\ref{f:config}, with an odd
number of cusps between the strands $U$ and $U'$. Then, $(Y_K,\xi_K)$
is overtwisted.
\begin{figure}[ht]
\begin{center}
\epsfig{file=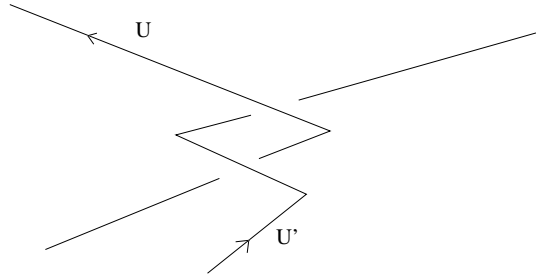, height=3.6cm}
\end{center}
\caption{Configuration producing an overtwisted disk}
\label{f:config}
\end{figure}
\end{thm}

\begin{cor}\label{c:negtorus}
Let $K$ be a Legendrian knot in the standard contact three--sphere. 
If $K$ is smoothly isotopic to a negative torus knot then 
$(Y_K, \xi_K)$ is overtwisted.
\end{cor}

Notice the contrast: when the Legendrian knot $K$ satisfies $\tb(K)=
2g_s(K)-1$ (where $g_s(K)$ denotes its slice genus) then $(Y_K,
\xi_K)$ is tight \cite{LS4}.  The tightness question for contact
structures can be fruitfully attacked with the use of the contact
Ozsv\'ath--Szab\'o invariants \cite{OSz6}.  In fact, the nonvanishing
of these invariants implies tightness, while their computation can
sometimes be performed (see e.g.~\cite{LS4,LS5}) using a contact
surgery presentation in conjunction with the surgery exact triangle
established in Heegaard Floer theory by Peter Ozsv\'ath and Zolt\'an
Szab\'o~\cite{OSzabs}. Such ideas can be used to prove the following.

\begin{thm}\label{t:ntb}
Let $K\subset S^3$ be a smooth knot. Suppose that, for some integer
$n>0$, the 3--manifold $S^3_n(K)$ is a lens space. Let $L\subset
(S^3,\xi_{st})$ be a Legendrian knot smoothly isotopic to $K$. Then,
$L$ has Thurston--Bennequin invariant not greater than $n$.
\end{thm}

In the proof of Theorem~\ref{t:ntb} we will only assume that $S^3_n
(K)$ is an $L$--space, a weaker condition specified in
Section~\ref{s:prelim} and known to be satisfied by lens spaces.

In our investigations we prove tightness by establishing the
nonvanishing of the appropriate contact Ozsv\'ath--Szab\'o invariant.
Therefore, we are interested in cases when this invariant vanishes,
although overtwistedness does not obviously hold.

\begin{prop}\label{p:vanish1}
Let $L_1, L_2\subset (S^3,\xi_{st})$ be two smoothly isotopic
Legendrian knots whose Thurston--Bennequin invariants satisfy
\[
\tb(L_1)<\tb(L_2).
\]
Then, the result of contact $(+1)$--surgery along $L_1$ has vanishing contact
Ozsv\'ath--Szab\'o  invariant.
If $\tb (L)\leq -2$ then the contact Ozsv\'ath--Szab\'o invariant 
$c^+(Y_L, \xi _L)$ vanishes.
\end{prop}

\begin{rem}
The hypotheses of Proposition~\ref{p:vanish1} do not imply that either
$L_1$ or $L_2$ be stabilizations of other Legendrian knots. In fact,
examples of Legendrian knots $L_1$ and $L_2$ satisfying the
assumptions of Proposition~\ref{p:vanish1} without being
stabilizations were found by Etnyre and Honda \cite{EHLA}.
\end{rem}

In many cases the contact invariants can be explicitly computed. We
will perform such computations for a subfamily of Legendrian knots
called Chekanov--Eliashberg knots, cf.~\cite{EFM}. These knots are of
particular interest because they have equal ``classical invariants''
(i.e., knot type, Thurston--Bennequin invariant and rotation number)
but are not Legendrian isotopic. Our computation shows that, at least
when combined with the particular surgery approach we adopt here, the
contact Ozsv\'ath--Szab\'o invariant is not strong enough to
distinguish these knots up to Legendrian isotopy.  For the precise
formulation of this fact see Section~\ref{s:checkanov}.

As a further application, we present examples where the contact
Ozsv\'ath--Szab\'o invariants distinguish contact structures defined
on a fixed 3--manifold.  In particular, by a simple calculation we
recover the main result of \cite{LM}:

\begin{thm}[\cite{LM}]\label{lm}
The Brieskorn integral homology sphere $-\Sigma (2,3,6n-1)$ admits
at least $(n-1)$ nonisotopic tight contact structures.
\end{thm}
\begin{rem} The same result was obtained in a more general form
by O.~Plamenevskaya \cite{OP}.
\end{rem}

Section~\ref{s:prelim} is devoted to the necessary (and brief)
recollection of background information about contact surgery and
Ozsv\'ath--Szab\'o invariants. Proofs of most of the statements
announced in the Introduction are given in
Section~\ref{s:proofs}. Section~\ref{s:checkanov} is devoted to 
the Legendrian Chekanov--Eliashberg
knots. In Section~\ref{s:dist} we prove Theorem~\ref{lm}.

\section{Preliminaries}\label{s:prelim}

For the basics of contact geometry and topology we refer the reader
to~\cite{Et,Ge}.

\subsection*{Contact surgery}
Let $(Y, \xi )$ be a closed, contact 3--manifold and $L\subset (Y,\xi
)$ a Legendrian knot. The contact structure $\xi$ can be extended from
the complement of a neighborhood of $L$ to the 3--manifold obtained by
$(\pm 1)$--surgery along $L$ (with respect to the contact framing). In
fact, by the classification of tight contact structures on the solid
torus $S^1\times D^2$ \cite{H1}, such an extension is uniquely
specified by requiring that its restriction to the surgered solid
torus be tight.  The same uniqueness property holds for all surgery
coefficients of the form $\frac{1}{k}$ with $k \in \Z$. For a general
nonzero rational surgery coefficient, there is a finite number of
choices for the extension.  Consequently, a Legendrian knot $L\subset
(S^3, \xi _{st})$ decorated with $+1$ or $-1$ gives rise to a
well--defined contact 3--manifold, which we shall denote by $(Y_L, \xi
_L)$ and $(Y^L, \xi ^L)$, respectively.  For a more extensive
discussion on contact surgery see \cite{DG2}.

\subsection*{Heegaard Floer theory}
In this subsection we recall the basics of the Ozsv\'ath--Szab\'o
homology groups.  For a more detailed treatment see \cite{OSzF1,
OSzF2, OSzF4}.

According~ \cite{OSzF1}, to a closed, oriented spin$^c$ 3--manifold
$(Y, \t )$ one can associate a finitely generated Abelian group $\hf
(Y, \t )$ and a finitely generated $\Z [U]$--module $HF^+(Y,
\t )$. A spin$^c$ cobordism $(W, \s)$ between $(Y_1,\t_1)$ and
$(Y_2,\t_2)$ gives rise to homomorphisms $\hat F_{W, \s}\colon \hf (Y_1, \t
_1)\to \hf (Y_2, \t _2)$ and $F^+_{W, \s }\colon HF^+(Y_1, \t _1)\to
HF^+(Y_2, \t _2)$, with $F^+_{W,\s}$ $U$--equivariant. 

Let $Y$ be a closed, oriented 3--manifold and $K\subset Y$ a
framed knot with framing $f$. Let $Y(K)$ denote the 3--manifold given
by surgery along $K\subset Y$ with respect to the framing $f$. The
surgery can be viewed at the 4--manifold level as a 2--handle
addition. The resulting cobordism $X$ induces a homomorphism
\[
\hat F_X:=\sum_{\s\in Spin^c(X)} \hat F_{X,\s}\co \hf (Y)\to\hf(Y(K)),
\]
where $\hf(Y):=\oplus_{\t\in Spin^c (Y)} \hf(Y,\t)$. 
Similarly, there is a cobordism $Z$ defined by adding a 2--handle to 
$Y(K)$ along a normal circle $N$ to $K$ 
with framing $-1$ with respect to a normal disk to $K$. 
The boundary components of $Z$ are $Y(K)$ and the
3--manifold $Y'(K)$ obtained from $Y$ by a surgery along $K$ with
framing $f+1$. As before, $Z$ induces a homomorphism
\[
\hat F_Z\co\hf(Y(K))\to\hf(Y'(K)).
\]
The above construction can be repeated starting with $Y(K)$ and
$N\subset Y(K)$ equipped with the framing specified above: we get $Z$
(playing the role previously played by $X$) and a new cobordism $W$
starting from $Y'(K)$, given by attaching a 4--dimensional 2--handle
along a normal circle $C$ to $N$ with framing $-1$ with respect to a
normal disk. It is easy to check that this last operation yields $Y$
at the 3--manifold level. 

\begin{thm} [\cite{OSzF2}, Theorem~9.16]\label{t:triangle}
The homomorphisms $\hat F_{X}, \hat F_Z$  and $\hat F_W$ fit 
into an exact triangle
\[
\begin{graph}(6,2)
\graphlinecolour{1}\grapharrowtype{2}
\textnode {A}(1,1.5){$\hf (Y)$}
\textnode {B}(5, 1.5){$\hf (Y(K))$}
\textnode {C}(3, 0){$\hf (Y'(K))$}
\diredge {A}{B}[\graphlinecolour{0}]
\diredge {B}{C}[\graphlinecolour{0}]
\diredge {C}{A}[\graphlinecolour{0}]
\freetext (3,1.8){$\hat F_X$}
\freetext (4.6,0.6){$\hat F_Z$}
\freetext (1.4,0.6){$\hat F_W$}
\end{graph}
\]
\qed
\end{thm}

For a torsion spin$^c$ structure (i.e.~a spin$^c$ structure whose
first Chern class is torsion) the homology theories $\hf$ and $HF^+$
come with a relative $\Z$--grading which admits a lift to an absolute
$\Q$--grading \cite{OSzabs}. The action of $U$ shifts this degree by
$-2$.


For $a\in \Q$, define $\T_a:=\oplus_b (\T_a)_b$ as the graded $\Z
[U]$--module such that, for every $b\in\Q$,
\[
(\T_a)_b=
\begin{cases}
\Z\quad\text{for}\quad b\geq a\quad\text{and}\quad b-a\in 2\Z,\\ 
0\quad\text{otherwise},
\end{cases}
\]
and the $U$--action $(\T _a)_b \to (\T _a)_{b-2}$ is an isomorphism
for every $b\neq a$. The following proposition can be extracted
from~\cite[Theorem~10.1]{OSzF2} and~\cite[Propositions~4.2
and~4.10]{OSzabs}.

\begin{prop}[\cite{OSzF2, OSzabs}]\label{p:struct}
Let $Y$ be a rational homology sphere. Then, for each $\t\in Spin^c (Y)$
\[
HF^+(Y, \t ) = \T _a \oplus A(Y),
\]
where $a\in \Q$ and $A(Y)=\oplus_d A_d(Y)$ is a graded, finitely
generated Abelian group. Moreover,
\[
HF^+(-Y,\t) = \T_{-a}\oplus A(-Y),
\]
with $A_d(-Y)\cong A_{-d-1}(Y)$. If $b_1 (Y)=1$ and $\t \in Spin^c (Y)$ is
torsion then
\[
HF^+(Y, \t )=\T _a \oplus \T _{a'} \oplus A'(Y), 
\]
where $a-a'$ is an odd integer and $A'(Y)=\oplus_d A'_d(Y)$ is a
graded, finitely generated Abelian group. Moreover,
\[
HF^+(-Y, \t )=\T _{-a} \oplus \T _{-a'} \oplus A'(Y), 
\]
with $A'_d(-Y)\cong A'_{-d-1}(Y)$.\qed
\end{prop}

The two theories $\hf$ and $HF^+$ are related by a long exact sequence, 
which takes the following form for a torsion spin$^c$ structure $\t$
\begin{equation}\label{e:exseq}
\ldots \to \hf _a(Y, \t ) \stackrel{f}{\longrightarrow} 
HF^+_a(Y, \t ) \stackrel{U}{\longrightarrow}  HF^+_{a-2}(Y, \t ) 
\to \hf _{a-1}(Y, \t) \to  \ldots
\end{equation}
where $U$ denotes ``multiplication by $U$''. All the gradings
appearing in the sequence can be worked out from the definitions and
the construction of the exact sequence (cf.~\cite[Section~2]{OSzabs}).

\begin{cor}\label{c:bgr}
Let $Y$ be a rational homology 3--sphere. Then, $HF^+(Y, \t )\cong\T
_a$ if and only if $\hf (Y, \t )\cong\Z $. If $b_1(Y)=1$ and $\t$ is a
torsion spin$^c$ structure, then $HF^+(Y, \t )\cong\T _{a_1}\oplus \T
_{a_2}$ if and only if $\hf (Y, \t )\cong\Z ^2$.
\end{cor}

\begin{proof}
We sketch the proof of the statement for $b_1(Y)=0$, the other case
can be proved by similar arguments. Clearly, if $HF^+(Y, \t )\cong\T _a$
then it follows immediately from Exact Sequence~\eqref{e:exseq} that $\hf
(Y, \t ) = \hf _a(Y, \t )\cong\Z$.  Conversely, if $\hf(Y, \t)\cong\Z$
then Exact Sequence~\eqref{e:exseq} and Proposition~\ref{p:struct} imply
$HF^+(Y, \t )\cong\T _a$.
\end{proof}

Observe that, in view of Corollary~\ref{c:bgr}, if $Y$ is a rational
homology 3--sphere, the following two conditions are equivalent:
\begin{enumerate}
\item
For each spin$^c$ structure $\t\in Spin ^c(Y)$, $HF^+(Y, \t )\cong\T
_a$ for some $a$;
\item
For each spin$^c$ structure $\t\in Spin ^c(Y)$, $\hf (Y, \t )\cong\Z$.
\end{enumerate}

\begin{defn}
A rational homology 3--sphere satisfying any of the above equivalent
conditions is called an \emph{$L$--space}.
\end{defn}

It follows from Proposition~\ref{p:struct} that an oriented rational
homology 3--sphere $Y$ is an $L$--space if and only if $-Y$ is an
$L$--space.  Moreover, lens spaces are
$L$--spaces~\cite[Section~3]{OSzF2}.

We will use the following fact regarding the maps connecting the
Ozsv\'ath--Szab\'o homology groups. Suppose that $W$ is a cobordism
defined by a single 2--handle attachment.

\begin{prop}[\cite{LS4}] \label{p:adjunction}
Let $W$ be a cobordism containing a smooth, closed, oriented surface
$\Sigma$ of genus $g$, with $\Sigma\cdot\Sigma > 2g-2$. Then, the
induced maps $\hat F_{W, \s }$ and $F^+_{W, \s }$ vanish for every spin$^c$
structures $\s$ on $W$. \qed
\end{prop}

\subsection*{Contact Ozsv\'ath--Szab\'o invariants}
Let $(Y,\xi)$ be a closed, contact 3--manifold. Then, the contact
Ozsv\'ath--Szab\'o invariants
\[
\chat (Y, \xi ) \in \hf (-Y, \t _{\xi })/\langle \pm
1\rangle\quad\text{and}\quad c^+(Y, \xi ) \in HF^+(-Y, \t_{\xi
})/\langle \pm 1\rangle
\]
are defined~\cite{OSz6}, with $f(\chat (Y, \xi ))=c^+(Y, \xi )$, where
$f$ is the homomorphism appearing in Exact Sequence~\eqref{e:exseq}
and $\t _{\xi}$ is the spin$^c$ structure induced by the contact
structure $\xi$.

To simplify notation, throughout the paper we ignore the sign
ambiguity in the definition of the contact invariants, and treat them
as honest elements of the appropriate homology groups rather than
equivalence classes. The reader should have no problem checking that
there is no loss in making this abuse of notation. Alternatively, one
could work with $\Z/2\Z$ coefficients to make the sign ambiguity
disappear altogether. The properties of $\chat$ and $c^+$ which will
be relevant for us can be summarized as follows.

\begin{thm}[\cite{OSz6}]\label{t:item}
Let $(Y,\xi)$ be a closed, contact 3--manifold, and denote by
$c(Y,\xi)$ either one of the contact invariants $\chat(Y,\xi)$ and
$c^+(Y,\xi)$. Then,
\begin{enumerate}
\item 
The class $c(Y, \xi )$ is an invariant of the isotopy class
of the contact structure $\xi $ on $Y$.
\item 
If $(Y, \xi )$ is overtwisted then $c(Y, \xi )=0$, while 
if $(Y, \xi )$ is Stein fillable then $c(Y, \xi )\neq 0$.
\item 
Suppose that $(Y_2,\xi_2)$ is obtained from $(Y_1,\xi _1)$ by a
contact $(+1)$--surgery. Then we have 
\[
F_{-X} (c(Y_1, \xi_1))= c(Y_2,\xi_2), 
\]
where $-X$ is the cobordism induced by the surgery with orientation
reversed and $F_{-X}$ is the sum of $F_{-X,\s}$ over all spin$^c$
structures $\s$ extending the spin$^c$ structures induced on $-Y_i$ by
$\xi_i$, $i=1,2$.  In particular, if $c(Y_2, \xi_2)\neq 0$ then $(Y_1,
\xi_1)$ is tight.
\item 
Suppose that $\t _{\xi }$ is torsion. Then
$c(Y, \xi )$ is a homogeneous element of degree
$-h(\xi )\in \Q$, where $h(\xi )$ is the Hopf--invariant of the
2--plane field defined by the contact structure $\xi $. \qed
\end{enumerate}
\end{thm}

\begin{rem}
The Hopf--invariant can be easily determined for a contact structure
defined by a contact $(\pm 1)$--surgery diagram along the Legendrian
link $\Li \subset (S^3, \xi _{st})$ \cite{DGS}. In fact, fix an
orientation of $\Li$ and consider the 4--manifold $X$ defined by the
Kirby diagram specified by the surgery \cite{GS}.  Let $c\in H^2 (X;
\bfz )$ denote the cohomology class which evaluates as rot$(L)$ on the
homology class determined by a component $L$ of the link $\Li$. If $\t
_{\xi }$ is torsion, then $c^2\in \Q$ is defined, and $h(\xi)$ is
equal to $\frac{1}{4}(c^2 -3\sigma (X) - 2\chi (X)+2)+q$, where $q$ is
the number of $(+1)$--surgeries made along $\Li$ to get $(Y, \xi )$.
\end{rem}

\section{Proofs}\label{s:proofs}
Now we can turn to the proofs of the statements announced in
Section~\ref{s:intro}.  

\begin{proof}[Proof of Theorem~\ref{t:ot}]
Consider the Legendrian push--off $K'$ of $K$ drawn as a dotted line
in the left--hand side of Figure~\ref{f:modification2}. 
\begin{figure}[ht]
\includegraphics[width=15cm]{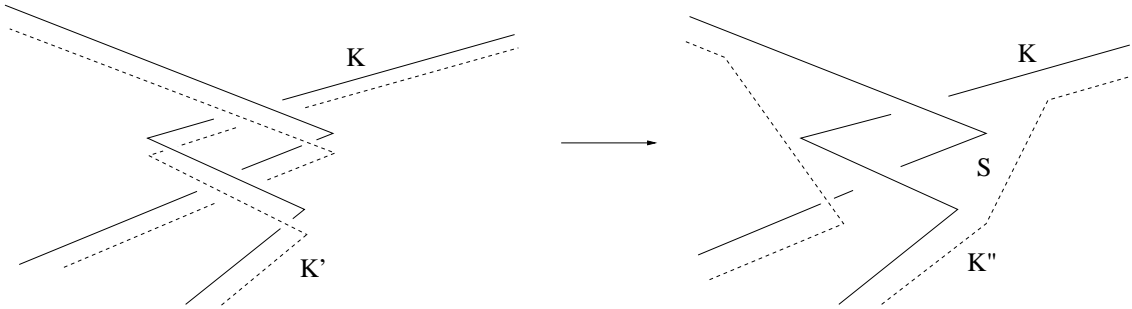}
\caption{The modification of the push--off}
\label{f:modification2} 
\end{figure}
The obvious annulus between $K$ and $K'$ induces framing $\tb(K)$ on
both $K$ and $K'$.  Consider the modification $K''$ of $K'$
illustrated in the right--hand side of
Figure~\ref{f:modification2}. The obvious surface $S$ between $K''$
and $K$ is oriented because of the hypotheses on the cusps of the
front projection, it has genus $1$ and it induces framing $\tb(K)+1$
on $K$. In particular, $S$ extends to a meridian disk $D$ inside the
surgered solid torus. Since $S$ induces framing $\tb(K)+1$ on $K''$,
while $\tb(K'')=\tb(K')+3=\tb(K)+3$, we have $\tb_{S\cup D}(K'')=2$,
i.e. the Legendrian knot $K''=\del (S\cup D)$ violates the
Bennequin--Eliashberg inequality with respect to $S\cup D$. We
conclude that $(Y_K, \xi _K)$ is overtwisted.
\end{proof} 

To prove Theorem~\ref{t:ntb}, Corollary~\ref{c:negtorus} and
Proposition~\ref{p:vanish1} we shall need the following lemma (for a
different proof of a more general result see~\cite{Oz}).

\begin{lem}\label{l:+1ot}
Let $K$ be a Legendrian knot in the standard contact three--sphere.
If $K$ is the stabilization of another Legendrian knot then
$(Y_K, \xi_K)$ is overtwisted. 
\end{lem}

\begin{proof}
By assumption, $K$ admits a front projection containing one of the
configurations of Figure~\ref{f:zig-zags}. Without loss, we may assume
that we are in the situation of the left--hand side of
Figure~\ref{f:zig-zags}.
\begin{figure}[ht]
\includegraphics[height=4cm]{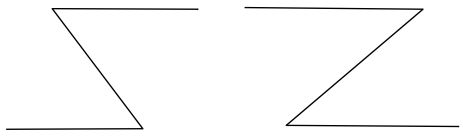}
\caption{The two possible ``zig--zags''}
\label{f:zig-zags} 
\end{figure}
Consider the Legendrian push--off $K'$ of $K$ drawn as a dotted line
in the left--hand side of Figure~\ref{f:modification1}. The obvious annulus
between $K$ and $K'$ induces framing $\tb(K)$ on both $K$ and $K'$.
Consider the modification $K''$ of $K'$ illustrated in the right--hand
side of Figure~\ref{f:modification1}. 
\begin{figure}[ht]
\includegraphics[height=4cm]{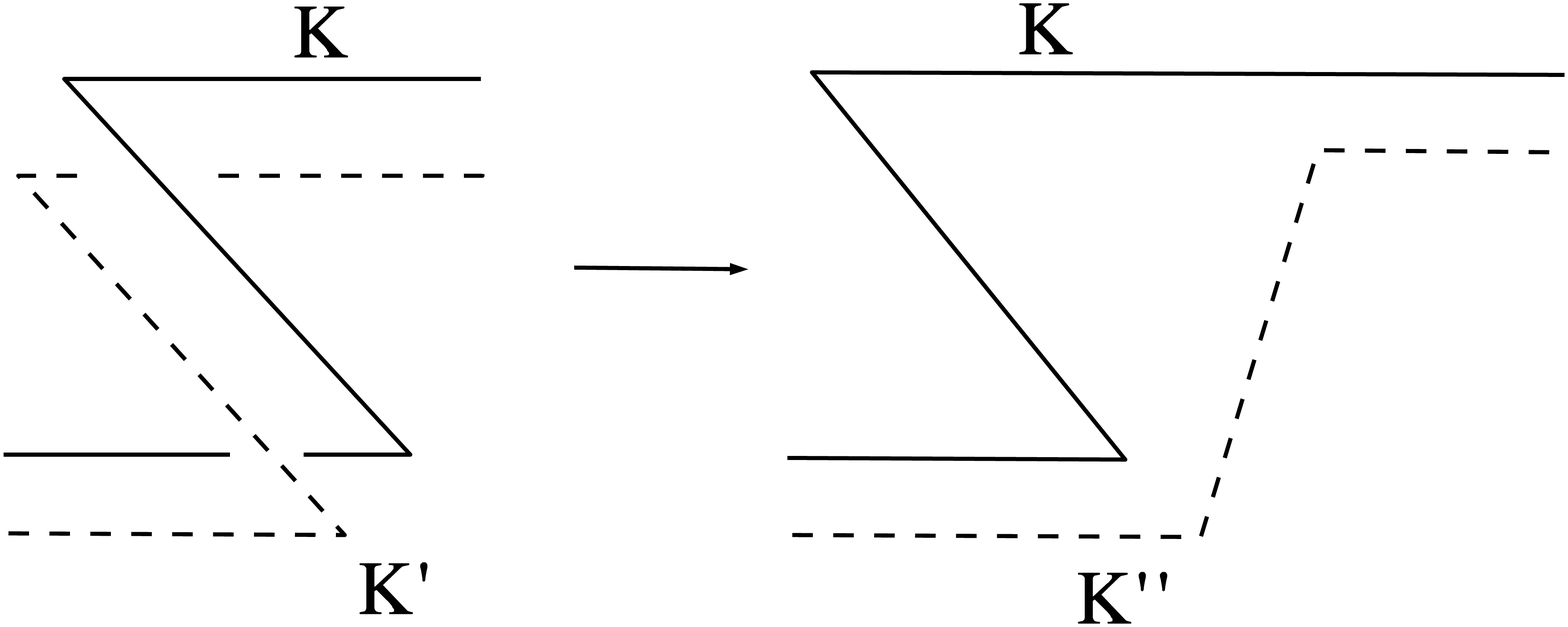}
\caption{The modification of the Legendrian push--off}
\label{f:modification1} 
\end{figure}
There still is an obvious annulus $A$ between $K''$ and $K$, except
that now it induces framing $\tb(K'')=\tb(K)+1$ on $K$ and $K''$.
Since we perform contact $(+1)$--surgery on $K$, the annulus $A$
extends to a meridian disk $D$ inside the surgered solid
torus. Therefore, $D\cup A$ is an overtwisted disk in $(Y_K, \xi
_K)$. 
\end{proof}

The proof of Lemma~\ref{l:+1ot} clearly applies to establish the
following slight generalization:

\begin{prop}
Suppose that the Legendrian link $\Li \subset (S^3, \xi _{st})$ is
obtained by stabilizing some components of another Legendrian
link. Let $(Y_{\Li }, \xi _{\Li} )$ be the result of contact $(\pm
1)$--surgeries along the components of $\Li$. If the surgery
coefficient on one of the stabilized components is $(+1)$, then
$(Y_{\Li }, \xi _{\Li })$ is overtwisted. \qed
\end{prop}

\begin{proof}[Proof of Corollary~\ref{c:negtorus}]
Examining~\cite[Figure~8]{EH3}, it is easy to check that any
Legendrian negative torus knot $K$ with maximal Thurston--Bennequin
invariant contains the configuration of Figure~\ref{f:config}, with an
odd number of cusps between the two strands $U$ and $U'$.  Therefore,
by Theorem~\ref{t:ot} $(Y_K,\xi_K)$ is overtwisted. On the other hand,
according to the results of~\cite{EH3}, any Legendrian negative torus
knot $K'$ with non--maximal Thurston--Bennequin invariant is isotopic
to the stabilization of one with maximal Thurston--Bennequin
invariant. Thus, by Lemma~\ref{l:+1ot} $(Y_{K'},\xi_{K'})$ is
overtwisted.
\end{proof}

\begin{proof}[Proof of Theorem~\ref{t:ntb}]
By contradiction, suppose that $S^3_n (K)$ is an $L$--space (recall
that lens spaces are $L$--spaces) and $L_1\subset (S^3, \xi _{st})$ is
a Legendrian knot smoothly isotopic to $K$ with $tb(L_1)>n$. Let $L$
be obtained by stabilizing $L_1$ $\tb(L_1)-n$ times, so that
$\tb(L)=n$. Denote by $(Y_L, \xi _L)$ the result of contact
$(+1)$--surgery along $L$. By Lemma~\ref{l:+1ot} $(Y_L, \xi _L)$ is
overtwisted, hence $\hat c(Y_L, \xi _L)=0$.  On the other hand, we can
compute $\hat c(Y_L, \xi _L)$ using Theorem~\ref{t:item}, getting
$\hat c(Y_L, \xi _L)=\hat F_{-X}(c(S^3, \xi _{st}))$ where $X$ is the
appropriate cobordism. The map $\hat F_{-X}$ fits into the exact
triangle
\[
\begin{graph}(6,2)
\graphlinecolour{1}\grapharrowtype{2}
\textnode {A}(1,1.5){$\hf (S^3)$}
\textnode {B}(5, 1.5){$\hf (S^3_{-n-1}(\mk ))$}
\textnode {C}(3, 0){$\hf (S^3 _{-n}(\mk ))$}
\diredge {A}{B}[\graphlinecolour{0}]
\diredge {B}{C}[\graphlinecolour{0}]
\diredge {C}{A}[\graphlinecolour{0}]
\freetext (2.9,1.9){$\hat F_{-X}$}
\freetext (1.4,0.6){$\hat F_W$}
\end{graph}
\]
where $\overline K$ is the mirror image of $K$ and $S^3_r(K)$ denotes
the result of $r$--surgery along $K$.  Since $S^3_{-n}(\mk )
=-S^3_n (K)$ is an
$L$--space, we have
\[
\rk\hf(S^3_{-n}(\overline K)) = |H_1(S^3_{-n}(\overline K))| = n,
\]
while by Proposition~\ref{p:struct}
\[
\rk\hf (S^3 _{-n-1}(\mk ))\geq 
|H_1(S^3 _{-n-1}(\mk ))|=n+1.
\]
Exactness of the triangle immediately implies $\hat F_W=0$, therefore
$\hat F_{-X}$ must be injective.  Since $\hat c (S^3, \xi _{st})\neq
0$, this shows $\hat c(Y_L, \xi _L)\neq 0$, which contradicts the fact
that $(Y_L, \xi _L)$ is overtwisted.
\end{proof}

\begin{proof}[Proof of Proposition~\ref{p:vanish1}] 
Consider a Legendrian knot $L'$ obtained by stabilizing $L_2$ until
$tb (L_1)=tb (L')$. Since $L'$ and $L_1$ are smoothly isotopic and
have the same contact framing, the cobordisms associated to the
contact $(+1)$--surgeries along $L_1$ and $L'$ can be
identified. Since $c(Y_{L_1}, \xi_{L_1})$ and $c(Y_{L'}, \xi _{L'})$
are images of $c(S^3, \xi _{st})$ under the same map, $c(Y_{L_1},
\xi_{L_1})=0$ if and only if $c(Y_{L'},
\xi_{L'})=0$. Lemma~\ref{l:+1ot} gives $c(Y_{L'}, \xi_{L'})=0$, and
the first statement follows.

For the second statement consider the exact triangle in the $HF^+$--theory
provided by the surgery along $L$. (The Thurston--Bennequin invariant 
$\tb (L)$ is denoted by $t$.) 
After reversing orientation the triangle takes the shape 
\[
\begin{graph}(6,2)
\graphlinecolour{1}\grapharrowtype{2}
\textnode {A}(1,1.5){$HF^+ (S^3)$}
\textnode {B}(5, 1.5){$HF^+ (S^3_{-t-1}(\ml ))$}
\textnode {C}(3, 0){$\hf (S^3 _{-t}(\ml ))$}
\diredge {A}{B}[\graphlinecolour{0}]
\diredge {B}{C}[\graphlinecolour{0}]
\diredge {C}{A}[\graphlinecolour{0}]
\freetext (2.9,1.9){$F^+_{-W}$}
\end{graph}
\]
Now the assumption $t<-1$ implies that $-t-1>0$, hence the cobordism
$-W$ inducing the first map is positive definite. 
It is known that the map $F^{\infty }_{-W}$ on the $HF^{\infty }$--theory
vanishes if $b_2^+(-W)>0$ \cite{OSzF2}. Since for $S^3$ the natural map
$HF ^{\infty }(S^3)\to HF^+(S^3)$ is onto, this implies that $F^+_{-W}=0$.
Since 
\[
c^+(Y_L, \xi _L)=F^{+}_{-W}(c^+(S^3, \xi _{st})),
\]
the vanishing of the contact invariant $c^+(Y_L, \xi _L)$ follows.
\end{proof}

\section{Examples}\label{s:checkanov}
Given a Legendrian knot $L\subset (S^3, \xi _{st})$, we shall denote
by $(Y_L, \xi _L)$, respectively $(Y^L, \xi ^L)$, the contact
3--manifold obtained by contact $(+1)$--, respectively
$(-1)$--surgery.

\begin{figure}[ht]
\begin{center}
\epsfig{file=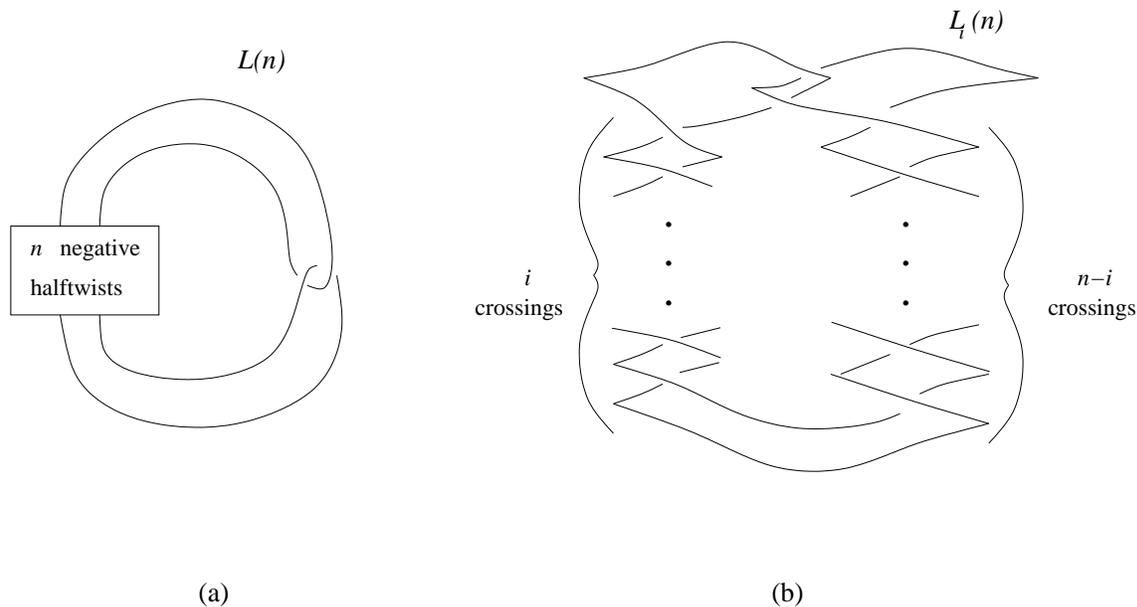, height=8cm}
\end{center}
\caption{The $n$-twist knot, and its Legendrian realizations}
\label{f:twist}
\end{figure}

Let $L_i=L_i (n)$, $i=1,\ldots,n-1$, be the Legendrian knot given by
Figure~\ref{f:twist}(b).  The knots $L_i(n)$ ($n$ fixed and $\geq 2$)
were considered in~\cite{EFM}. They are all smoothly isotopic to the
$n$-twist knot of Figure~\ref{f:twist}(a) (having $n$ negative
half--twists).
The knots $L_i$ were the first examples of smoothly isotopic
Legendrian knots having equal classical invariants
(i.e. Thurston--Bennequin invariants and rotation numbers) but not
Legendrian isotopic \cite{Check, EFM}.  The reader should be aware
that our convention for representing a Legendrian knot via its front
projection differs from the one used in~\cite{EFM}. In fact, we use
the contact structure given by the $1$--form $dz+xdy$ rather than the
$1$--form $-dz+ydx$, used in~\cite{EFM}.  However, the
contactomorphism between the two contact structures given by sending
$(x,y,z)$ to $(y,-x,z)$ induces a one--to--one correspondence between
the corresponding front projections, and under this correspondence
Figure~1 from~\cite{EFM} is sent to our Figure~\ref{f:twist}(b).

\begin{prop}\label{p:check1}
For every $1\leq i, j \leq n-1$ we have 
\[
\hat c(Y_{L_i}, \xi _{L_i}) = \hat c(Y_{L_j}, \xi _{L_j}).
\]
\end{prop}

\begin{proof}
The statement follows easily from basic properties of the contact
invariant: by the surgery formula for contact $(+1)$--surgeries, we
have $\hat c(Y_{L_i}, \xi _{L_i})= F_{-X}(\hat c(S^3,\xi _{st}))$,
where $X$ is the cobordism induced by the 4--dimensional handle
attachment dictated by the surgery. Since $X$ depends only on the
smooth isotopy class of the Legendrian knot and its
Thurston--Bennequin invariant, and is therefore independent of $i$,
the claim trivially follows.
\end{proof}

According to the main result of this section, Theorem~\ref{t:check2},
the same equality holds if we perform Legendrian surgeries along
$L_i(n)$, that is, the contact Ozsv\'ath--Szab\'o invariants of the
results of contact $(\pm 1)$--surgeries do not distinguish the
Chekanov--Eliashberg knots.

\begin{thm}\label{t:check2}
Let $n\geq 2$ be an even integer, and let $1\leq i,j\leq n-1$ be both
odd. Then,
\[
\hat c(Y^{L_i}, \xi ^{L_i}) = \hat c(Y^{L_j}, \xi ^{L_j}).
\]
\end{thm}

The proof of Theorem~\ref{t:check2} rests on the following 
two lemmas.  

\begin{lem}[\cite{OSzalt}] \label{l:alt}
Let $n\geq 2$ be an even integer, and denote by $\overline L(n)$ the
mirror image of $L(n)$. Then,
\[
HF^+(S_0^3({\overline {L}}(n)))\cong\T _{\frac{1}{2}}\oplus \T _{\frac{3}{2}}
\oplus \Z ^{\frac n2 - 1}_{(\frac{1}{2})}. 
\]
\end{lem}

\begin{proof}
Let $k=\frac n2$. Choosing a suitable oriented basis for an obvious
Seifert surface for $L(n)$ one can easily compute the Seifert matrix
\[
\begin{pmatrix}
-k & k-1\\
k & -k
\end{pmatrix}, 
\]
with eigenvalues $-1$ and $1-4k$. This immediately gives signature
$\sigma (L(n))=-2$ and Alexander polynomial
\[
\Delta _{L(n)}(t) = k t^{-1}-(2k-1) + k t.
\]
Since $L(n)$ is an alternating knot with genus $g(L(n))=1$,
applying~\cite[Theorem~1.4]{OSzalt} we get
\[
\begin{cases}
HF^+(S^3 _0 (L(n)), \s)\cong\T _{-\frac{1}{2}}
\oplus \T _{-\frac{3}{2}}\oplus \Z ^{\frac n2 - 1}_{(-\frac{3}{2})}
\quad\text{if}\quad c_1(\s)=0,\\
HF^+(S^3 _0 (L(n)), \s )=0\quad\text{if}\quad c_1(\s )\neq 0.
\end{cases}
\]
By Proposition~\ref{p:struct} this implies the result.
\end{proof}

\begin{figure}[ht]
\begin{center}
\epsfig{file=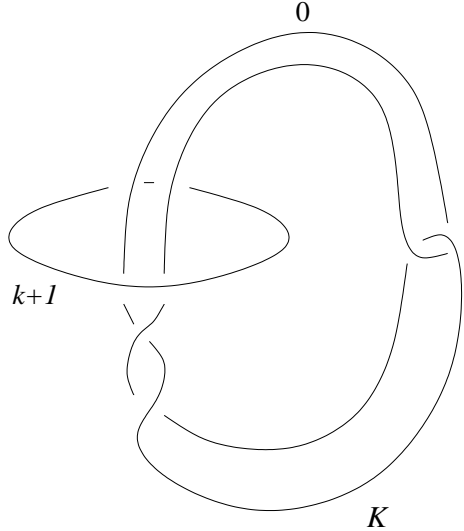, height=7cm}
\end{center}
\caption{Surgery diagram for $V(k)$}
\label{f:kirby}
\end{figure}

\begin{lem} \label{l:comp}
Let $k\geq 0$ be an integer, and let $V(k)$ be the oriented 3--manifold defined by the surgery diagram of Figure~\ref{f:kirby}. Then,  
\[
\hf (V(k)) \cong \Z ^{2k+2}\quad\text{and}\quad
HF^+(V(k))=\oplus _{i=1}^{2k+2} \T _{a_i}\quad
\text{for some $a_i \in \Q$}.
\]
\end{lem}

\begin{proof}
In order to compute $\hf (V(k))$ we will use the exact triangle
defined by the $(k+1)$--framed unknot of Figure~\ref{f:kirby}. 
It is easy to see that the unknot of Figure~\ref{f:kirby} bounds 
a punctured torus smoothly embedded in the complement of 
the knot $K$. Thus, the cobordism we get by
attaching this last 2--handle contains a torus with self--intersection
$(k+1)$, and the induced map in the surgery triangle vanishes by
Proposition~\ref{p:adjunction}. Consequently, the surgery triangle is
actually a short exact sequence.  Notice that $K$ is the
(left--handed) trefoil knot, hence $\hf (S^3_0(K))=\Z
^2$~\cite[Theorem~1.4]{OSzalt}. Arguing by induction we get
\[
\hf (V(k+1))\cong\hf (V(k))\oplus \Z ^2
\]
for every $k\geq 0$. On the other hand, for $k=0$ the unknot can be blown 
down,
showing that $V(0)\cong S^1\times S^2$. This fact immediately implies 
\begin{equation}\label{e:hhatvn}
\hf (V(k))\cong\Z ^{2k+2}
\end{equation}
for every $k\geq 0$.
Using the surgery presentation of Figure~\ref{f:kirby} it is easy to
check that
\[
H_1(V(k); \Z )\cong \Z \oplus \Z/(k+1)\Z,
\]
therefore $V(k)$ admits $(k+1)$ different torsion spin$^c$
structures. By Proposition~\ref{p:struct} and Exact Sequence~\eqref{e:exseq} 
we have 
\[
\rk\hf(V(k),\t)\geq 2
\]
if $\t$ is a torsion spin$^c$ structure. 
Therefore, using~\eqref{e:hhatvn}, we see that $\hf (V(k), \t )\cong\Z^2$ for
each torsion spin$^c$ structure $\t$ and  
\[
\hf(V(k),\t)=0
\]
if $\t$ is not torsion.  The statement now follows from
Proposition~\ref{p:struct} and Corollary~\ref{c:bgr}.
\end{proof}

\begin{proof}[Proof of Theorem~\ref{t:check2}]
The idea of the proof is the following: First we will find a contact
3--manifold $(Y, \xi )$ such that contact $(+1)$--surgery along some
Legendrian knot $K\subset (Y, \xi )$ gives $(Y^{L_i}, \xi ^{L_i})$ and
$A(Y)\subset HF^+(Y, \t _{\xi })$ (as it is defined in
Proposition~\ref{p:struct}) vanishes. Therefore $c^+(Y, \xi )$ is an
element of some $\T _a$. The $U$--equivariance of the map induced by
the surgery will then show that $c^+(Y^{L_i}, \xi ^{L_i}) \in \T _a
\subset HF^+(Y^{L_i}, \t _{\xi ^{L_i}})$, from which the conclusion
will easily follow.

To this end, consider the contact structure $\eta _i (n)$ defined by
Legendrian surgery along the 2--component link of Figure~\ref{f:link}.
\begin{figure}[ht]
\begin{center}
\epsfig{file=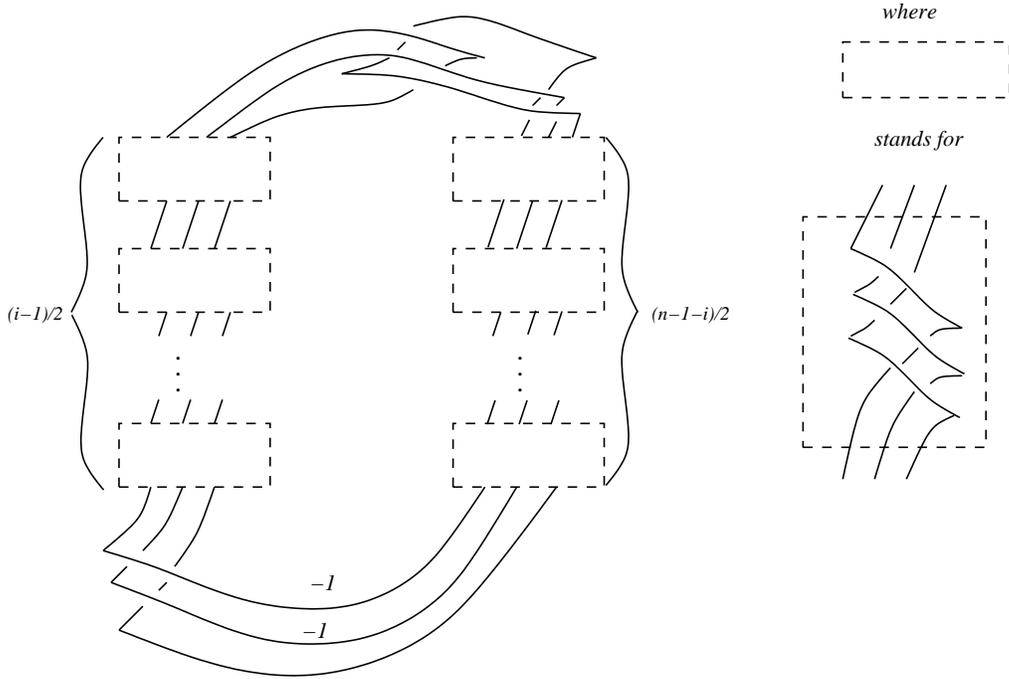, height=9cm}
\end{center}
\caption{Contact surgery diagram defining $(Y(n), \eta _i (n))$}
\label{f:link}
\end{figure}
Notice that one of the knots in Figure~\ref{f:link} is topologically
the unknot, while the other one is $L_i(n)$.  According to the Kirby
moves indicated in Figure~\ref{f:moves}, it follows that this contact
structure lives on the 3--manifold $Y(n):=-V(\frac n2)$, where $V(k)$
is defined by
\begin{figure}[ht]
\begin{center}
\epsfig{file=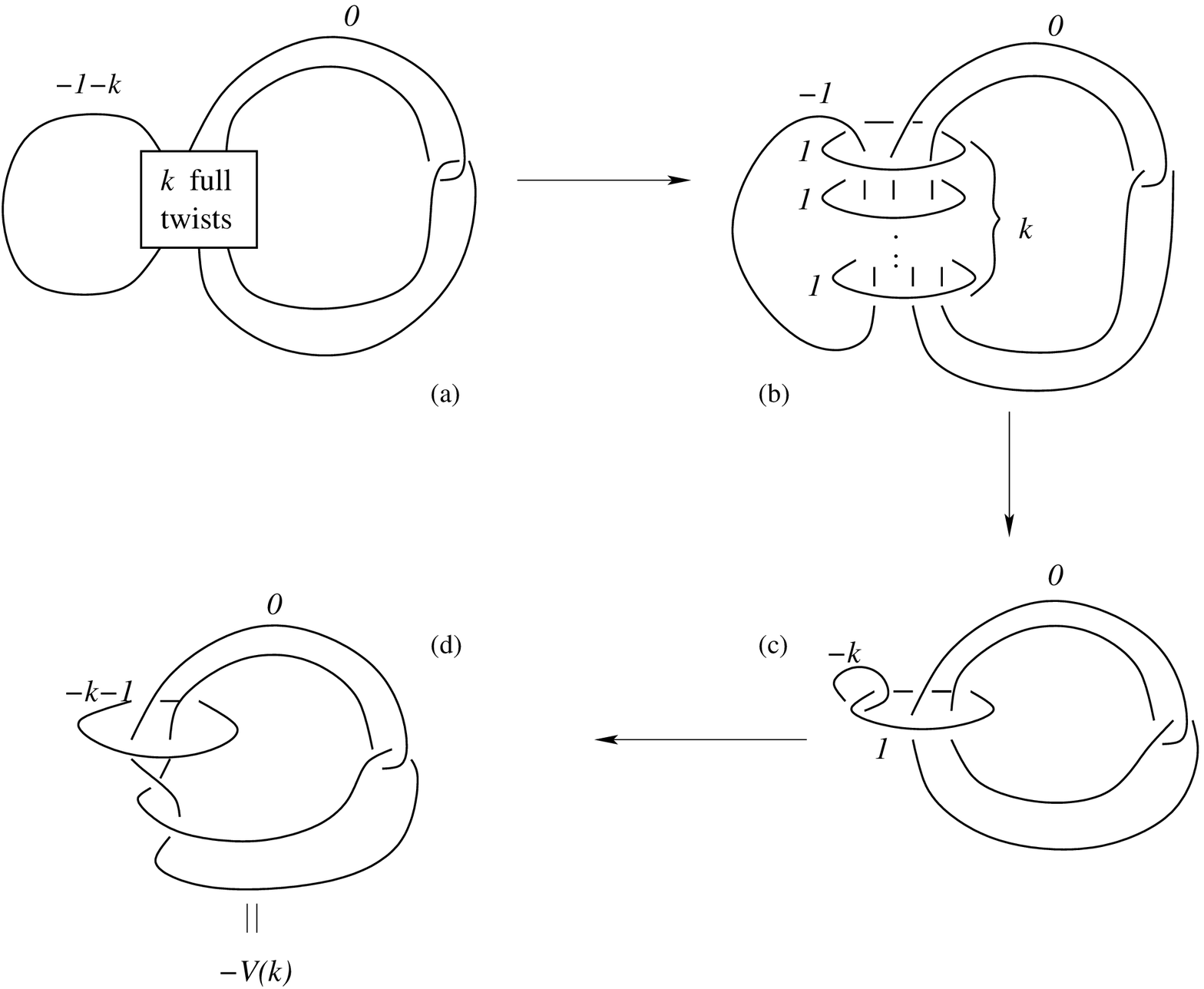, height=11cm}
\end{center}
\caption{Kirby moves for $Y(n)$}
\label{f:moves}
\end{figure}
Figure~\ref{f:kirby}. According to \cite{DG1}, the effect of a contact
$(\pm 1)$--surgery along a Legendrian knot can be cancelled by contact
$(\mp 1)$--surgery along a Legendrian push--off of the knot.
Therefore, doing contact $(+1)$--surgery along the push--off of the
unknot in Figure~\ref{f:link} we get $(Y^{L_i}, \xi ^{L_i})$.  On the
other hand, denoting by $X_n$ the cobordism induced by the contact
$(+1)$-surgery, we have
\[
\hat F_{-X_n}(\hat c(Y(n),\eta _i(n)))=\hat c(Y^{L_i}, \xi ^{L_i}).
\]
A simple computation shows that $h(\xi^{L_i})=-\frac{1}{2}$, 
therefore by Theorem~\ref{t:item}(4) we have 
\[
\hat c(Y^{L_i}, \xi ^{L_i})\in\hf_{\frac{1}{2}}( -Y^{L_i}). 
\]
Moreover, $\hat c(Y^{L_i}, \xi ^{L_i})$ is primitive~\cite{OP}.  Thus,
to prove the statement it will be enough to verify that there is a
rank--$1$ subgroup of $\hf_{\frac{1}{2}}( -Y^{L_i})$ containing $\hat
F_{-X_n}(\hat c(Y(n),\eta _i(n)))$ for every $i$.  An easy computation
shows that (since we assumed $n$ to be even)
the Thurston--Bennequin numbers of the knots $L_i(n)$ are
all equal to 1, cf. \cite{EFM}, hence each of the 3--manifolds $Y^{L_i}$ is
diffeomorphic to $S^3_0(L(n))$. By Lemma~\ref{l:alt}
\[
HF^+(-S_0 ^3 (L (n)))\cong\T _{\frac{1}{2}} \oplus \T _{\frac{3}{2}} \oplus A, 
\]
where $A$ is a finitely generated abelian group, while by Lemma~\ref{l:comp} we have
\[
HF^+(-Y(n))=\oplus_{i=1}^{n+2} \T _{a_i}
\]
for some $a_i\in \Q$. Since $F^+_{-X_n}$ is
$U$--equivariant and for sufficiently large $h$ the action of $U^h$
vanishes on $A$, we have
\[
{\rm Im} (F^+_{-X_n}) \subseteq 
\T _{\frac{1}{2}} \oplus \T _{\frac{3}{2}} \subseteq HF^+(-S_0 ^3 (L(n))).  
\]
Therefore, up to sign, there is a unique primitive element in ${\rm
Im}(F^+_{-X_n})$ of degree $\frac{1}{2}$, implying that $c^+(Y^{L_i},
\xi ^{L_i})=c^+ (Y^{L_j}, \xi ^{L_j})$ for $i,j$ as in the statement.
Since
\[
HF^+_{-\frac{1}{2}}(-S_0^3(L(n)))=0, 
\]
it follows that the homomorphism 
\[
f\co\hf _{\frac{1}{2}}(-S_0^3(L(n)))\to 
HF^+_{\frac{1}{2}}(-S_0^3(L(n)))
\]
from Exact Sequence~\eqref{e:exseq} is injective. Since 
\[
f(\hat c(Y^{L_i}, \xi ^{L_i}))=c^+(Y^{L_i}, \xi ^{L_i})\in {\rm Im}(F^+_{-X_n})
\]
for every $i$, this concludes the proof. 
\end{proof}

\section{Distinguishing tight contact structures}\label{s:dist}

\begin{defn}
Let $\xi_i$, for $i=1,\ldots, n-1$, denote the contact structure on the
Brieskorn sphere $-\Sigma (2,3,6n-1)$ defined by the contact surgery specified by Figure~\ref{f:structures}.
\begin{figure}[ht]
\includegraphics[height=7cm]{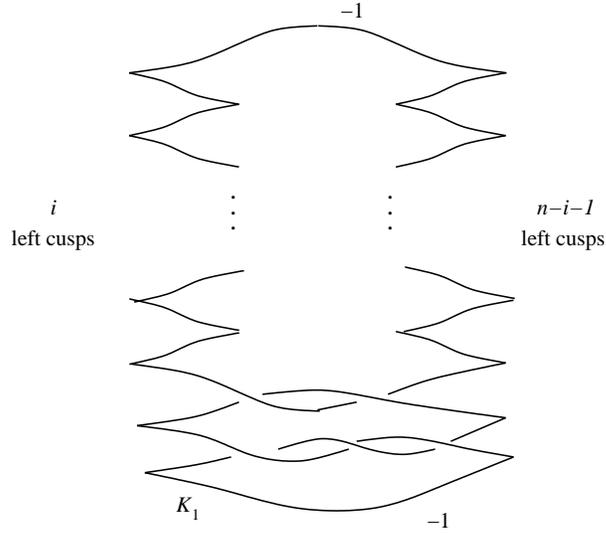}
\caption{Contact structures on the 3--manifold $-\Sigma (2,3,6n-1)$}
\label{f:structures}
\end{figure}
\end{defn}

\begin{thm}\label{t:lm1}
The contact invariants $c^+(\xi_1), \ldots , c^+(\xi_{n-1})$
are linearly independent over $\Z$.
\end{thm}

\begin{proof}
Consider the Legendrian push-off ${\tilde {K}}_1$ of the Legendrian
trefoil $K_1$ of Figure~\ref{f:structures}.  Attach a 4--dimensional
2--handle along ${\tilde {K}}_1$ to $-\Sigma (2,3,6n-1)$ with framing
equal to the contact framing $+1$.  Since contact $(+1)$--surgery
along a Legendrian push--off cancels contact $(-1)$--surgery, we get a
cobordism $W$ such that $F_{-W}(c^+(\xi_i))= c^+(\eta _i)$, where
$\eta_i $ is the contact structure on $L(n, 1)$ defined by
Figure~\ref{f:lens}.
\begin{figure}[ht]
\includegraphics[height=8cm]{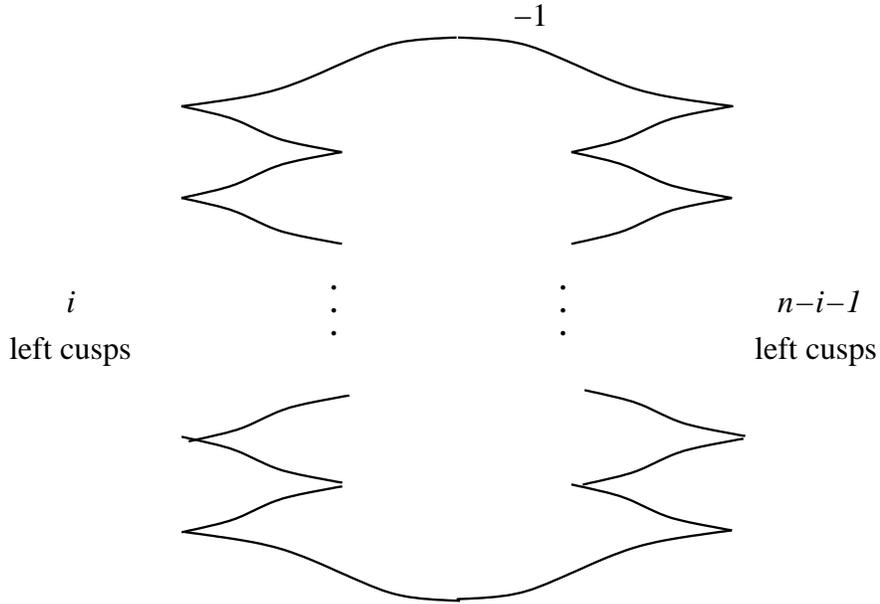}
\caption{The contact structure $\eta_i$ on $L(n,1)$}
\label{f:lens}
\end{figure}
The contact invariants $c^+(\eta _i)$ are linearly independent
because they belong to groups corresponding to different spin$^c$
structures on the same lens space $L(n,1)$. Therefore, the invariants
$c^+(\xi_i)$ are also linearly independent, concluding the proof.
\end{proof}

\begin{cor}\label{c:lm1}
The contact structures $\xi_1,\ldots,\xi_{n-1}$ are pairwise non--isotopic.
\qed\end{cor}

Corollary~\ref{c:lm1} was first proved by Lisca and Mati\'c \cite{LM}
using Seiberg--Witten theory. For a different Heegaard Floer theoretic
proof (of a more general statement) see \cite{OP}.

\begin{rem}
It is known \cite{OSzabs} that $HF^+(-\Sigma (2,3,6n-1)) =\T _{-2}
\oplus \Z _{(-2)} ^{n-1}$, therefore by Proposition~\ref{p:struct}
$HF^+ (\Sigma (2,3,6n-1))=\T _2 \oplus \Z _{(1)} ^{n-1}$.  It follows
from Theorem~\ref{t:lm1} that the elements $c^+(\xi _i)$ ($i=1, \ldots
, n-1$) span $HF^+_1(\Sigma (2,3,6n-1))$.
\end{rem}

Notice that if the trefoil knot of Figure~\ref{f:structures} is
replaced by any Legendrian knot $L$, the statement of
Theorem~\ref{t:lm1} holds with the same proof. If tb$(L)=1$ and
rot$(L)=0$, then the contact resulting structures $\xi_1, \ldots , \xi
_{n-1}$ are all homotopic as 2--plane fields.

\end{document}